\documentclass{article}
\usepackage{amsmath}
\usepackage{amsfonts}
\usepackage{amssymb}

\newtheorem{theorem}{Theorem}

\newtheorem{corollary}[theorem]{Corollary}

\newtheorem{lemma}[theorem]{Lemma}

\newenvironment{proof}[1][Proof]{\textbf{#1.} }{\ \rule{0.5em}{0.5em}}
\input{tcilatex}

\begin{document}

\title{On a class of linearizable planar geodesic webs}
\author{Vladislav V. Goldberg and Valentin V. Lychagin}
\date{}
\maketitle

\begin{abstract}
We present a complete description of a class of linearizable planar geodesic
webs which contain a parallelizable $3$-subweb.
\end{abstract}

\section{Introduction}

The paper is a continuation of \cite{GL09b}. In the paper \cite{GL09b} we
considered some classical problems of the theory of planar webs. In
particular, at the end of the paper we proved that \emph{a planar $d$-web is
linearizable if and only if the web is geodesic and the Liouville tensor of
one of its $4$-subwebs vanishes.} In the current paper we describe all
linearizable planar geodesic webs satisfying the following additional
condition: the curvature $K$ of one of its 3-subwebs vanishes.

\section{The Problem}

Below we give some (not all) definitions and notions which will be used in
the paper. For additional information a reader is advised to look into \cite%
{GL09b}.

We consider the plane $M$ endowed with a torsion-free connection $\nabla $
and a geodesic $d$-web in $M$, i.e., a $d$-web all leaves of all foliations
of which are geodesic with respect to the connection $\nabla $. We have
proved in \cite{GL09b} that \emph{there is a unique projective structure
associated with a planar $4$-web in such a way that the $4$-web is geodesic
with respect to the structure.}

The flatness of the projective structure can be checked by the Liouville
tensor (see \cite{Lio89}, \cite{Lie83}, \cite{K08}).  This tensor can be
constructed as follows (see, for example, \cite{NS94}).

Let $\nabla$ be a representative of the canonical projective structure, and $%
Ric$ be the Ricci tensor of the connection $\nabla$. Define a new tensor $%
\mathfrak{P}$ as
\begin{equation*}
\mathfrak{P}(X,Y)=\frac{2}{3}Ric(X,Y)+\frac{1}{3}Ric(Y,X),
\end{equation*}
where $X$ and $Y$ are arbitrary vector fields.

The Liouville tensor $\mathfrak{L}$ is defined as follows:
\begin{equation*}
\mathfrak{L}(X,Y,Z)=\nabla_{X}(\mathfrak{P})(Y,Z)-\nabla_{Y}(\mathfrak{P}%
)(X,Z)
\end{equation*}
where $X, Y$ and $Z$ are arbitrary vector fields.

The tensor is skew-symmetric in $X$ and $Y$, and therefore it belongs to
\begin{equation*}
\mathfrak{L}\in\Omega^{1}(\mathbb{R}^{2})\otimes\Omega^{2}(\mathbb{R}^{2}).
\end{equation*}
It is known (see \cite{Lio89}, \cite{NS94}, \cite{Lie83}, \cite{K08}) that
\emph{the Liouville tensor depends on the projective structure defined by $%
\nabla$ and vanishes if and only if the projective structure is flat.}

For the case of the projective structure associated with a planar 4-web we
shall call this tensor the \emph{Liouville tensor} of the 4-web.

Let us consider a $4$-web with a $3$-subweb given by a web function $f\left(
x,y\right) $ and a basic invariant $a$ (see \cite{GL09b} for more details) and
introduce the following three invariants:
\begin{equation}
w=\frac{f_{y}}{f_{x}},\;\;\alpha =\frac{aa_{y}-wa_{x}}{wa(1-a)},\;\;k=(\log
w)_{xy}.  \label{invar}
\end{equation}

Then the Liouville tensor has the form \cite{GL09b}:
\begin{equation*}
\mathfrak{L}=(L_{1}\omega _{1}+\frac{L_{2}}{w}\omega _{2})\otimes \omega
_{1}\wedge \omega _{2},
\end{equation*}%
where $L_{1}$ and $L_{2}$ are relative differential invariants of order
three.

The explicit formulas for these invariants are
\begin{equation}
\renewcommand{\arraystretch}{1.5}%
\begin{array}{lll}
3L_{1} & = & w(-(kw)_{x}+\alpha _{xx}+\alpha \alpha _{x})+(\alpha
w_{xx}+(\alpha ^{2}+3\alpha _{x})w_{x}-2\alpha _{xy}-2\alpha \alpha _{y}) \\
&  & +w^{-1}(-\alpha w_{xy}-2\alpha _{y}w_{x}+\alpha w_{x}^{2})+w^{-2}\alpha
w_{x}w_{y}, \\
3L_{2} & = & w^{2}(-(kw^{-1})_{y}+2\alpha \alpha _{x})+w(2\alpha
^{2}w_{x}-2\alpha _{xy}-\alpha \alpha _{y}) \\
&  & +(-\alpha w_{xy}-2\alpha _{y}w_{x}+\alpha _{yy})+w^{-1}(\alpha
w_{x}w_{y}-\alpha _{y}w_{y}).%
\end{array}%
\renewcommand{\arraystretch}{1}  \label{Lio tensor}
\end{equation}

As we said in Introduction, at the end of the paper \cite{GL09b} we proved
that a planar $d$-web is linearizable if and only if the web is geodesic and
the Liouville tensor of one of its $4$-subwebs vanishes.

In the current paper we consider \emph{a class of planar $d$-webs for which
the curvature $K$ of one of its $3$-subwebs vanishes.}

In order to prove the main theorem, we need the following lemma.

\begin{lemma}
If $K = 0$, we can reduce $w$ $($see $(\ref{invar}))$ to one: $w=1$.
\end{lemma}

\begin{proof}
In fact, because
\begin{equation*}
K=-\displaystyle\frac{1}{f_{x}f_{y}}\Biggl(\log \displaystyle\frac{f_{x}}{%
f_{y}}\Biggr)_{xy},
\end{equation*}%
it follows from $K=0$ that $(\log w)_{xy}=0$. Hence $\log w=u(x)+v(y)$,
where $u(x)$ and $v(y)$ are arbitrary functions. It follows that $w=a(x)b(y)$, where $a(x)=e^{u(x)}$ and $b(y)=e^{v(y)}$. Taking the gauge
transformation $x\rightarrow X\left(x\right), y\rightarrow Y\left(y\right)
$, with $X^{\prime }\left( x\right) =$ $e^{u(x)}$ and $Y^{\prime}\left(y\right)
=e^{-v(y)}$, we get that $w=1$.
\end{proof}

We shall prove now the main theorem.

\begin{theorem}
A planar $d$-web, for which the curvature $K$ of one of its $3$-subwebs
vanishes, is linearizable if and only if the web is geodesic,  and the
invariants $\alpha $ defined  by its $4$-subwebs have one of the following
forms:

\begin{description}
\item[$(i)$]
\begin{equation}
\alpha =\frac{\wp ^{\prime }(2x+y+\lambda _{1},g_{2},g_{3})-\wp ^{\prime
}(x+2y+\lambda _{2},g_{2},-g_{3})}{\wp (2x+y+\lambda _{1},g_{2},g_{3})-\wp
(x+2y+\lambda _{2},g_{2},-g_{3})},  \label{type1}
\end{equation}%
where $\wp $ is the Weierstrass function, $g_{2}$ and $g_{3}$ are
invariants, and $\lambda _{1}$ and $\lambda _{2}$ are arbitrary constants.

\item[$(ii)$]
\begin{equation}
\alpha =k\frac{e^{k(x-y+C)}+1}{e^{k(x-y+C)}-1},  \label{type2}
\end{equation}%
where $k$ and $C$ are arbitrary constants.

\item[$(iii)$]
\begin{equation}
\alpha =-k\tan \frac{x-y+C}{2},  \label{type3}
\end{equation}%
where $k$ and $C$ are arbitrary constants.

\item[$(iv)$]
\begin{equation}
\alpha =\frac{2}{x-y+C},  \label{type4}
\end{equation}%
where $C$ is an arbitrary constant.
\end{description}

Here $x,y$ are such coordinates that the $3$-subweb is defined by
the web functions $x, y$ and $x+y.$
\end{theorem}

\begin{proof}
By Theorem 9 of \cite{GL09b}, the conditions of linearizability are $%
L_{1}=0,L_{2}=0$. By (\ref{invar}) and Lemma 1, the condition $K=0$ implies $%
k=0,w=1$.

It follows that the conditions $L_{1}=0,L_{2}=0$ become
\begin{equation}
\renewcommand{\arraystretch}{1.5}\left\{
\begin{array}{ll}
\alpha _{xx}-2\alpha _{xy}+\alpha \alpha _{x}-2\alpha \alpha _{y}=0, &  \\
\alpha _{yy}-2\alpha _{xy}+2\alpha \alpha _{x}-\alpha \alpha _{y}=0. &
\end{array}%
\right. \renewcommand{\arraystretch}{1}  \label{L=0}
\end{equation}%
Conditions (\ref{L=0}) can be written in the form
\begin{equation}
\renewcommand{\arraystretch}{1.5}\left\{
\begin{array}{ll}
(\partial _{x}-2\partial _{y})(\alpha _{x}+\frac{1}{2}\alpha ^{2})=0, &  \\
(\partial _{y}-2\partial _{x})(\alpha _{y}-\frac{1}{2}\alpha ^{2})=0. &
\end{array}%
\right. \renewcommand{\arraystretch}{1}  \label{L=0/2}
\end{equation}%
Therefore, relations (\ref{L=0/2}) imply
\begin{equation}
\renewcommand{\arraystretch}{1.5}\left\{
\begin{array}{ll}
\alpha _{x}+\frac{1}{2}\alpha ^{2}=A(2x+y), &  \\
\alpha _{y}-\frac{1}{2}\alpha ^{2}=B(x+2y) &
\end{array}%
\right. \renewcommand{\arraystretch}{1}  \label{A,B}
\end{equation}%
for some functions $A$ and $B$.

Differentiating the first equation of (\ref{A,B}) with respect to $y$ and
the second one with respect to $x$, we get the following compatibility
conditions for (\ref{A,B}):
\begin{equation*}
\alpha\alpha_{y}+\alpha\alpha_{x}=A^{\prime }-B^{\prime },
\end{equation*}
which by (\ref{A,B}) is equivalent to
\begin{equation}
(A+B)\alpha=A^{\prime }-B^{\prime }.  \label{alpha}
\end{equation}

We assume that $A+B \neq 0$. (The case $A+B=0$ will be considered
separately.) Then equation (\ref{alpha}) implies
\begin{equation}
\alpha = \frac{A^{\prime }-B^{\prime }}{A+B}.  \label{alpha1}
\end{equation}

Next, we substitute $\alpha$ from (\ref{alpha1}) into equations (\ref{L=0}).
As a result, we obtain that
\begin{equation}
\renewcommand{\arraystretch}{1.5} \left\{
\begin{array}{ll}
(2A^{\prime \prime }-B^{\prime \prime })(A+B)-(A^{\prime }-B^{\prime
})(2A^{\prime }+B^{\prime })+\frac{1}{2}(A^{\prime }-B^{\prime 2} &  \\
= A(A+B)^2, &  \\
(A^{\prime \prime }-2B^{\prime \prime })(A+B)-(A^{\prime }-B^{\prime
})(A^{\prime }+2B^{\prime })-\frac{1}{2}(A^{\prime }-B^{\prime 2} &  \\
= B(A+B)^2. &
\end{array}
\right. \renewcommand{\arraystretch}{1}  \label{L=0/3}
\end{equation}

Adding and subtracting equations (\ref{L=0/3}), we find that
\begin{equation}
(A^{\prime \prime }-B^{\prime \prime })(A+B)-(A^{\prime 2}-B^{\prime 2})=%
\frac{(A+B)^{3}}{3}  \label{L=0/4}
\end{equation}%
and
\begin{equation}
A^{\prime \prime }+B^{\prime \prime }=A^{2}-B^{2}.  \label{L=0/5}
\end{equation}

Therefore,
\begin{equation}
\renewcommand{\arraystretch}{1.5}\left\{
\begin{array}{l}
A^{\prime \prime}-A^{2}=c,   \\
(B^{\prime \prime }+B^{2}=-c,
\end{array}%
\right. \renewcommand{\arraystretch}{1}  \label{L=0/6}
\end{equation}%
for a constant $c\in \mathbb{R}$ .

Multiplying equations (\ref{L=0/6}) by $A^{\prime }$ and $B^{\prime }$,
respectively, we get
\begin{eqnarray*}
\renewcommand{\arraystretch}{1.3}
A^{\prime }A^{\prime \prime }-A^{\prime 2} =cA^{\prime}; \\
B^{\prime }B^{\prime \prime }+B^{\prime 2} =-cB^{\prime},
\renewcommand{\arraystretch}{1}
\end{eqnarray*}%
and
\begin{eqnarray*}
\renewcommand{\arraystretch}{1.5}
\Bigl(\frac{1}{2}A^{\prime 2}-\frac{1}{3}A^{\prime 3}-cA\Bigr)^{\prime}
=0; \\
\Bigl(\frac{1}{2}B^{\prime 2}+\frac{1}{3}B^{\prime 3}+cB\Bigr)^{\prime}
=0,
\renewcommand{\arraystretch}{1}
\end{eqnarray*}%
respectively.

This means that
\begin{equation}
\frac{1}{2}A^{\prime 2}-\frac{1}{3}A^{\prime 3}-cA=a(s)  \label{a(s)}
\end{equation}%
and
\begin{equation}
\frac{1}{2}B^{\prime 2}+\frac{1}{3}B^{\prime 3}+cB=b(t),  \label{b(t)}
\end{equation}%
where $s=x+2y$ and $t=2x+y$.

Now equations (\ref{a(s)}), (\ref{b(t)}) and (\ref{L=0/4}) give
\begin{equation}
b=a=\operatorname{const}. \in \mathbb{R}  \label{a=b}
\end{equation}

Remind that  solutions of the equation
\begin{equation}
y^{\prime 2}=4y^{3}-g_{2}y-g_{3}  \label{W}
\end{equation}%
have the form
\begin{equation}
y=\wp (x+\lambda ,g_{2},g_{3}),  \label{Wsol}
\end{equation}%
where $\wp $ is the Weierstrass function, $g_{2}$ and $g_{3}$ are
invariants, and $\lambda $ is an arbitrary constant.

By (\ref{a=b}), equations (\ref{a(s)}) and (\ref{b(t)}) can be written as
\begin{equation}
\renewcommand{\arraystretch}{1.5}
\begin{array}{ll}
A^{\prime 2}= \displaystyle\frac{2}{3} A^3 + 2cA + 2a,   \\
B^{\prime 2}= -\displaystyle\frac{2}{3} B^3 - 2cB - 2a.   \\
\end{array}
\renewcommand{\arraystretch}{1}  \label{W1}
\end{equation}
Taking $A = \beta \wp$ and $B = \gamma \wp$, substituting them into (\ref{W1}%
) and comparing the result with (\ref{W}), we find that
\begin{equation*}
\beta = 6, \gamma = -6; g_2 = -\frac{c}{3}, g_3 = -\frac{a}{18},
\end{equation*}
i.e., $g_2$ and $g_3$ are the same for both equations (\ref{W1}).

By (\ref{Wsol}), the solutions of (\ref{W1}) are
\begin{equation}
\renewcommand{\arraystretch}{1.5}\left\{
\begin{array}{ll}
A=6\wp (t+\lambda _{1},g_{2},g_{3}), &  \\
B=-6\wp (t+\lambda _{2},g_{2},-g_{3}), &
\end{array}%
\right. \renewcommand{\arraystretch}{1}  \label{W2}
\end{equation}%
where $g_{2}$ and $g_{3}$ are arbitrary constants.

Equations (\ref{W2}) can be now written as
\begin{equation}
\renewcommand{\arraystretch}{1.5}\left\{
\begin{array}{ll}
A=6\wp (2x+y+\lambda _{1},g_{2},g_{3}), &  \\
B=-6\wp (x+2y+\lambda _{2},g_{2},-g_{3}). &
\end{array}%
\right. \renewcommand{\arraystretch}{1}  \label{W3}
\end{equation}

Finally, equations (\ref{alpha1}) and (\ref{W3}) give the following
expression (\ref{type1}) for the invariant $\alpha $:
\begin{equation*}
\alpha =\frac{\wp ^{\prime }(2x+y+\lambda _{1},g_{2},g_{3})-\wp ^{\prime
}(x+2y+\lambda _{2},g_{2},-g_{3})}{\wp (2x+y+\lambda _{1},g_{2},g_{3})-\wp
(x+2y+\lambda _{2},g_{2},-g_{3})}.
\end{equation*}

Consider now the cases for which $A + B = 0$, i.e., the cases
\begin{equation*}
A = v, \;\; B = - v, \;\; v \in \mathbb{R}.
\end{equation*}

Then system (\ref{A,B}) has the form
\begin{equation}
\renewcommand{\arraystretch}{1.5} \left\{
\begin{array}{ll}
\alpha_{x} +\frac{1}{2}\alpha^2=v, &  \\
\alpha_{y} -\frac{1}{2}\alpha^2=-v &
\end{array}
\right. \renewcommand{\arraystretch}{1}  \label{v,-v}
\end{equation}
and is consistent.

It follows from (\ref{v,-v}) that  $\alpha_x + \alpha_y = 0$. The solution
of this equation is  $\alpha = \alpha (x - y)$. As a result, we can write
two equations  (\ref{v,-v}) as one equation
\begin{equation}
\alpha^{\prime }+\frac{1}{2}\alpha^2 = v.  \label{v}
\end{equation}

Three cases are possible:

\begin{description}
\item[$(ii)$] $v = \frac{1}{2} k^2, \; k \neq 0$. Then the solution of (\ref%
{v}) has the form (\ref{type2}).

\item[$(iii)$] $v = -\frac{1}{2} k^2, \; k \neq 0$. Then the solution of (%
\ref{v}) has the form (\ref{type3}).

\item[$(ii)$] $v = 0$. Then the solution of (\ref{v}) has the form (\ref%
{type4}).
\end{description}
\end{proof}

\begin{corollary}
If for a geodesic $d$ web the  basic invariants are solutions of the Euler
equation and one of its $3$-subwebs is parallelizable, then this web is linearizable.
\end{corollary}


\begin{thebibliography}{9}
\bibitem{GL07} Goldberg,~V.~V., Lychagin,~V.~V., \emph{Abelian equations and
rank problems for planar webs} (Russian), Izv. Vyssh. Uchebn. Zaved. Mat.
\textbf{2007}, no. 10, 40--76. English translation in Russian Math. (Iz.
VUZ) \textbf{51}, no. 11, 39--75 (2007). MR{\footnote{{\footnotesize In
the bibliography we will use the following abbreviations for the
review journals: JFM for Jahrbuch f\"{u}r die Fortschritte der
Mathematik, MR for Mathematical Reviews, and Zbl for Zentralblatt
f\"{u}r Mathematik.}}}2381928
(2008k:53029)

\bibitem{GL08a} Goldberg,~V.~V., Lychagin,~V.~V., \emph{Geodesic webs on a
two-dimensional manifold and Euler equations}, Acta Appl. Math.  (2009) (to
appear); see also arXiv: 0810.5392, pp. 1--15 (2009)

\bibitem{GL09b} Goldberg,~V.~V., Lychagin,~V.~V., \emph{On rank problems for
planar webs and projective structures}, in \emph{The Abel Symposium} 2008,
Springer (2009) (to appear); see also arXiv: 0812.0125v2, pp. 1--31 (2009)

\bibitem{K08} Kruglikov,~B., \emph{Point classification of $2$nd order ODEs:
Tresse classification revisited and beyond}, arXiv: 0809.4653, pp. 1--22
(2008)

\bibitem{Lie83} Lie, S., \emph{Klassification und Integration von gew\"{o}%
hnlichen Differentialgleichungen zwischen $x,y$, die eine Gruppe von
Transformation gestatten. III}, Archiv f\"{u}r Math. og Naturvidenskab
\textbf{8} (Kristiania, 1883), 371--458 (JFM \textbf{15}, p. 751); see also Gesammelte Abhandlungen.
Bd. 5 (1924), paper XIV, 362--427. JFM \textbf{50}, p. 2.

\bibitem{Lio89} Liouville,~R., \emph{M\'{e}moire sur les invariants de certaines \'{e}%
quations diff\'{e}rentielles et sur leurs applications}, Journal de l'\'{E}%
cole Polytechnique \textbf{59}, 7--76 (1889).
JFM \textbf{21}, p. 317

\bibitem{NS94} Nomizu,~K., Sasaki,~T., \emph{Affine Differential Geometry},
Cambridge Tracts in Mathematics, \textbf{111}, Cambridge University Press,
Cambridge (1994). MR1311248 (96e:53014); Zbl 834:53002
\end{thebibliography}
\end{document}